\documentclass[aos,MSNbibl,nameyear,seceqn,dvips]{arximspdf}
\usepackage{dcolumn}
\usepackage{graphicx}

\doi{10.1214/13-AOS1143} 
\volume{41}
\issue{4}
\pubyear{2013}
\firstpage{2176}
\lastpage{2196}

\makeatletter

\newcolumntype{d}[1]{D{.}{.}{#1}}

\newcommand{\cal}{\mathcal}

\newtheorem{theorem}{Theorem}[section]
\newtheorem{lemma}{Lemma}[section]
\newtheorem{cor}{Corollary}[section]

\makeatother

\begin{document}
\begin{frontmatter}

\title{Empirical likelihood on the full parameter space}
\runtitle{Extended empirical likelihood}

\begin{aug}
\author[A]{\fnms{Min} \snm{Tsao}\corref{}\thanksref{t1}\ead[label=e1]{mtsao@uvic.ca}}
\and
\author[A]{\fnms{Fan} \snm{Wu}\ead[label=e2]{fwu@uvic.ca}}
\runauthor{M. Tsao and F. Wu}
\affiliation{University of Victoria}
\address[A]{Department of Mathematics\\
\quad and Statistics\\
University of Victoria\\
Victoria, British Columbia\\
Canada V8W 3R4\\
\printead{e1}\\
\hphantom{E-mail: }\printead*{e2}} 
\end{aug}

\thankstext{t1}{Supported by a research grant from the
National Science and Engineering Research Council of Canada.}

\received{\smonth{2} \syear{2013}}
\revised{\smonth{4} \syear{2013}}

%
\begin{abstract}
We extend the empirical likelihood of Owen [\textit{Ann. Statist.}
\textbf{18} (1990) 90--120] by partitioning its domain into the
collection of its contours and \mbox{mapping} the contours through a
continuous sequence of similarity transformations onto the full
parameter space. The resulting extended empirical likelihood is a
natural generalization of the original empirical likelihood to the full
parameter space; it has the same asymptotic properties and identically
shaped contours as the original empirical likelihood. It can also
attain the second order accuracy of the Bartlett corrected empirical
likelihood of DiCiccio, Hall and Romano [\textit{Ann. Statist.}
\textbf{19} (1991) 1053--1061]. A simple first order extended
empirical likelihood is found to be substantially more accurate than
the original empirical likelihood. It is also more accurate than
available second order empirical likelihood methods in most small
sample situations and competitive in accuracy in large sample
situations. Importantly, in many one-dimensional applications this
first order extended empirical likelihood is accurate for sample sizes
as small as ten, making it a practical and reliable choice for small
sample empirical likelihood inference.
\end{abstract}

%

\end{frontmatter}

\section{Introduction}

Since the seminal work of Owen (\citeyear{Owe88}, \citeyear{Owe90}),
there have been many advances in empirical likelihood method that have
brought applications of the method to virtually every area of
statistical research. It has been widely observed [e.g.,
\citet{HalLaS90,QinLaw94,CorDavSpa95,Owe01,LiuChe10}] that
empirical likelihood ratio confidence regions can have poor accuracy,
especially in small sample and multidimensional situations. In
particular, there is a persistent undercoverage problem in that
coverage probabilities of empirical likelihood ratio confidence regions
tend to be lower than the nominal levels. In this paper, we tackle a
fundamental problem underlying the poor accuracy and undercoverage,
that is, the empirical likelihood is defined on only a part of the
parameter space. We call this the \emph{mismatch problem} between the
domain of the empirical likelihood and the parameter space. We solve
this problem through a geometric approach that expands the domain to
the full parameter space. Our solution brings about substantial
improvements in accuracy of the empirical likelihood inference and is
particularly useful for small sample and multidimensional situations.

To see the mismatch problem, consider the example of empirical
likelihood for the mean based on $n$ observations $X_1, \ldots, X_n$ of
a random vector $X\in\mathbb{R}^d$. The underlying parameter space
$\Theta$ is $\mathbb{R}^d$ itself. But for a given $\theta\in\mathbb
{R}^d$, when the convex hull of the $(X_i-\theta)$ does not contain
$0$, the empirical likelihood $L(\theta)$, the empirical likelihood
ratio $R(\theta)=n^nL(\theta)$ and the empirical log-likelihood ratio
$l(\theta)=-2\log R(\theta)$ are all undefined.
When this occurs, an established convention assigns $L(\theta)=0$, and
technically $L(\theta)$, $R(\theta)$ and $l(\theta)$ are now defined
over the full parameter space. Nevertheless, to highlight the
difference between the natural domain of the empirical likelihood and
the parameter space, we define the common domain of $L(\theta)$,
$R(\theta)$ and $l(\theta)$ as
%
\begin{equation}\label{1-2}
\Theta_n=\bigl\{\theta\dvtx  \theta\in\Theta\mbox{ and } l(\theta)<+
\infty\bigr\}.
\end{equation}
With this definition, the mismatch can now be expressed as $\Theta
_n\subset\Theta$. For the mean, $\Theta_n$ is the interior of the
convex hull of the $X_i$, which is indeed a proper subset of $\Theta
=\mathbb{R}^d$. The mismatch $\Theta_n\subset\Theta$ holds for
empirical likelihoods in general as the basic formulation common to all
empirical likelihoods has a convex hull constraint on the origin, such
as the one for the mean above, which may be violated by some $\theta$
values in the parameter space $\Theta$.

The convex hull constraint violation underlying the mismatch is well
known in the empirical likelihood literature. It was first noted in
\citet{Owe90} for the case of the mean. See also
\citet{Owe01}. To assess its impact on coverage probabilities of
empirical likelihood ratio confidence regions, \citet{Tsa04}
investigated bounds on coverage probabilities resulting from the convex
hull constraint. To bypass this constraint, \citet{Bar07}
introduced a \emph{penalized empirical likelihood} (PEL) for the mean
which removes the convex hull constraint in the formulation of the
original empirical likelihood (OEL) of Owen (\citeyear{Owe90},
\citeyear{Owe01}) and replaces it with a penalizing term based on the
Mahalanobis distance. For parameters defined by general estimating
equations, \citet{CheVarAbr08} introduced an \emph{adjusted
empirical likelihood} (AEL) which retains the formulation of the OEL
but adds a pseudo-observation to the sample. The AEL is just the OEL
defined on the augmented sample. But due to the clever construction of
the pseudo-observation, the convex hull constraint will never be
violated by the AEL. \citet{EmeOwe09} showed that the AEL
statistic has a boundedness problem which may lead to trivial 100\%
confidence regions. They proposed an extension of the AEL involving
adding two pseudo-observations to the sample to address the boundedness
problem. \citet{CheHua13} also addressed the boundedness problem
by modifying the adjustment factor in the pseudo-observation.
\citet{LiuChe10} proved a surprising result that under a certain
level of adjustment, the AEL confidence region achieves the second
order accuracy of the Bartlett corrected empirical likelihood (BEL)
region by \citet{DiCHalRom91N2}. Recently, \citet{LahMuk12}
showed that under certain dependence structures, a~modified PEL for the
mean works in the extremely difficult case of large dimension and small
sample size. The PEL and the AEL are both defined on $\mathbb{R}^d$,
and are thus free from the mismatch problem.

In this paper, we propose a new \emph{extended empirical likelihood}
(EEL) that is also free from the mismatch problem. We derive this EEL
through the domain expansion method of \citet{Tsa13} which expands the
domain of the OEL but retains its important geometric characteristics.
This EEL makes effective use of the dimension information in the data
and can attain the second order accuracy of the BEL. The most important
aspect of this EEL, however, is that there is an easy-to-use first
order version which is substantially more accurate than the OEL. This
first order EEL is also more accurate than available second order
empirical likelihood methods in most small sample and multidimensional
applications and is comparable in accuracy to the latter when the
sample size is large. The focus of the present paper is on the
construction of EEL for the mean through which we introduce the basic
idea of and important tools for expanding the OLE domain to the full
parameter space. Under certain conditions, EEL for other parameters may
also be constructed but this will be discussed elsewhere.

For brevity, we will use ``OEL $l(\theta)$'' and ``EEL $l^*(\theta)$''
to refer to the original and extended empirical log-likelihood ratios
for the mean, respectively. Throughout this paper, we assume that the
parameter space $\Theta$ is $\mathbb{R}^d$. The case where $\Theta$ is
a known subset of $\mathbb{R}^d$ can be handled by finding EEL
$l^*(\theta)$ defined on $\mathbb{R}^d$ first and then, for $\theta
\notin\Theta$ only, redefine it as $l^*(\theta)=+\infty$.

\section{Preliminaries} We review several key results and assumptions
for developing the EEL defined on the full parameter space.

\subsection{Empirical likelihood for the mean} 

Let $X\in\mathbb{R}^d$ be a random vector with mean $\theta_0$ and
covariance matrix $\Sigma_0$. Two assumptions we will need are:
\begin{longlist}[$(A_2)$]
\item[$(A_1)$] $\Sigma_0$ is a finite covariance matrix with full rank
$d$; and

\item[$(A_2)$] $\lim\sup_{\|t\|\rightarrow\infty}|E[\exp\{it^TX\}]|
<1$ and $E\|X\|^{15}<+\infty$.
\end{longlist}
These are also assumptions under which the OEL for the mean is Bartlett
correctable [\citet{DiCHalRom91N1}, \citet{r4}].

Let $X_1, \ldots, X_n$ be independent copies of $X$ where $n>d$. Let
$\Theta_n$ be the collection of points in the interior of the convex
hull of the $X_i$. For a $\theta\in\mathbb{R}^d$, \citet{Owe90} defined
the empirical likelihood ratio $R(\theta)$ as
%
\begin{equation}\label{2-50}
R(\theta) = \sup\Biggl\{ \prod_{i=1}^{n}
nw_i \Big| \sum_{i=1}^nw_i(X_i-
\theta)=0, w_i \geq0, \sum_{i=1}^nw_i=1
\Biggr\}.
\end{equation}
It may be verified that $0<R(\theta)\leq1$ iff $\theta\in\Theta_n$.
Also, $R(\theta)=0$ if $\theta\notin\Theta_n$. Hence, the domain of
the OEL $l(\theta)=-2\log R(\theta)$ is $\Theta_n$.
For a $\theta\in\Theta_n$, the method of Lagrange multipliers may be
used to show that
%
\begin{equation}\label{2-58}
l(\theta) = 2 \sum^n_{i=1} \log\bigl\{1+
\lambda^T(X_i-\theta)\bigr\},
\end{equation}
where the multiplier $\lambda=\lambda(\theta)\in\mathbb{R}^d$ satisfies
%
\begin{equation}\label{2-56}
\sum_{i=1}^n \frac{X_i-\theta}{1+\lambda^T(X_i-\theta)} = 0.
\end{equation}
Under assumption ($A_1$), \citet{Owe90} showed that OEL $l(\theta)$ satisfies
%
\begin{equation}\label{2-60}
l(\theta_0) \stackrel{\mathcal{D}} {\longrightarrow}
\chi^2_d\qquad\mbox{as $n\rightarrow+\infty$}.
\end{equation}
For an $\alpha\in(0,1)$, let $c$ be the ($1-\alpha$)th quantile of
the $\chi^2_d$ distribution. Then, the $100(1-\alpha)\%$ OEL confidence
region for $\theta_0$ is given by
%
\begin{equation}\label{1-3}
{\mathcal C}_{1-\alpha}=\bigl\{\theta\dvtx  \theta\in\Theta_n \mbox{
and } l(\theta) \leq c\bigr\}.
\end{equation}
Under assumptions ($A_1$) and ($A_2$),
DiCiccio, Hall and Romano (\citeyear{DiCHalRom91N1,DiCHalRom91N2}) showed that the coverage error
of ${\mathcal C}_{1-\alpha}$ is $O(n^{-1})$, that is,
%
\begin{equation}\label{2-64}
P(\theta_0 \in{\mathcal C}_{1-\alpha})=P\bigl(l(
\theta_0)\leq c\bigr)=P\bigl(\chi^2_d \leq c
\bigr)+O\bigl(n^{-1}\bigr).
\end{equation}
More importantly, they showed that the empirical likelihood is Bartlett
correctable. To give a brief account of this surprising result, let
%
\begin{equation}\label{2-68}
{\mathcal C}_{1-\alpha}'=\bigl\{\theta\dvtx  l(\theta)
\bigl(1-bn^{-1}\bigr) \leq c\bigr\}
\end{equation}
be the Bartlett corrected empirical likelihood ratio confidence region
where $b$ is the \emph{Bartlett correction constant} and $(1-bn^{-1})$
is the \emph{Bartlett correction factor}, DiCiccio, Hall and Romano
(\citeyear{DiCHalRom91N1,DiCHalRom91N2}) showed that ${\mathcal
C}_{1-\alpha}'$ has a coverage error of only $O(n^{-2})$, that is,
%
\begin{equation}\label{2-70}
P\bigl(\theta_0 \in{\mathcal C}_{1-\alpha}'
\bigr)=P\bigl[l(\theta_0) \bigl(1-bn^{-1}\bigr)\leq c
\bigr]=P\bigl(\chi^2_d \leq c\bigr)+O\bigl(n^{-2}
\bigr).
\end{equation}
In practice, the Bartlett correction constant $b$ cannot be determined
since it depends on the moments of $X$ which are not available in the
nonparametric setting of the empirical likelihood. However, replacing
the Bartlett correction factor in (\ref{2-70}) with
$[1-bn^{-1}+O_p(n^{-3/2})]$ does not affect the $O(n^{-2})$ term in its
right-hand side, that is,
%
\begin{equation}\label{2-72}
P\bigl\{l(\theta_0)\bigl[1-bn^{-1}+O_p
\bigl(n^{-3/2}\bigr)\bigr]\leq c\bigr\}=P\bigl(\chi^2_d
\leq c\bigr)+O\bigl(n^{-2}\bigr).
\end{equation}
This allows us to replace $b$ in (\ref{2-68}) and (\ref{2-70}) with a
$\sqrt{n}$-consistent estimate $\hat{b}$ without invalidating (\ref
{2-70}). See \citet{DiCHalRom91N2} and \citet{HalLaS90} for detailed discussions on Bartlett correction.

\subsection{Extended empirical likelihood}

The OEL confidence region ${\mathcal C}_{1-\alpha}$ in (\ref{1-3}) is
confined to the OEL domain $\Theta_n$. This is a main cause of the
undercoverage problem associated with ${\mathcal C}_{1-\alpha}$ [\citet{Tsa04}]. To alleviate the problem, \citet{Tsa13} proposed to expand $\Theta
_n$ which will lead to larger EL confidence regions. Let $h_n\dvtx  \mathbb
{R}^d \rightarrow\mathbb{R}^d$ be a bijective mapping and define a
new empirical log-likelihood ratio $l^*(\theta)$ through the OEL
$l(\theta)$ as follows:
%
\begin{equation}\label{1-5}
l^*(\theta)=l\bigl(h^{-1}_n(\theta)\bigr)\qquad
\mbox{for $\theta\in\mathbb{R}^d$}.
\end{equation}
Then the domain for the new empirical log-likelihood ratio is $\Theta
_n^*= h_n(\Theta_n)$.
Here, $h_n$ plays the role of reassigning or extending the OEL values
of points in $\Theta_n$ to points in $\Theta_n^*$. Because of this,
\citet{Tsa13} named $l^*(\theta)$ the extended empirical log-likelihood
ratio or simply EEL. In particular, \citet{Tsa13} used the following
$\tilde{\theta}$-centred similarity mapping $h_n^*\dvtx \mathbb
{R}^d\rightarrow\mathbb{R}^d$
%
\begin{equation}\label{2-80}
h_n^*(\theta)=\tilde{\theta}+\gamma_n(\theta-\tilde{
\theta}),
\end{equation}
where $\tilde{\theta}$ is the sample mean and $\gamma_n\in\mathbb
{R}^1$ is a constant (which we will refer to as the \emph{expansion
factor}) satisfying $\gamma_n \geq0$ and $\gamma_n \rightarrow1$ as
$n \rightarrow+\infty$. If we choose $\gamma_n>1$, then $\Theta_n
\subset\Theta^*_n\subset\mathbb{R}^d$, and $\Theta^*_n$ alleviates the
mismatch problem of $\Theta_n$. The EEL confidence region for $\theta
_0$ is given by
%
\begin{equation}\label{1-7}
{\mathcal C}_{1-\alpha}^*=\bigl\{\theta\dvtx  \theta\in\Theta^*_n
\mbox{ and } l^*(\theta) \leq c\bigr\}.
\end{equation}

The advantages of the EEL based on $h_n^*$ in (\ref{2-80}) are: (1) the
EEL confidence regions are similarly transformed OEL confidence
regions, as such they retain the natural centre and shape of the OEL
confidence regions, (2) the EEL can be applied to empirical likelihood
inference for a wide range of parameters, and (3) with a properly
selected constant $\gamma_n$, EEL confidence regions can achieve the
second order accuracy of $O(n^{-2})$.

Nevertheless, the EEL based on $h_n^*$ is only a partial solution to
the mismatch problem because the domain of this EEL $\Theta^*_n$ is
also a proper subset of $\mathbb{R}^d$. A~second order version of this
EEL has been found to have good accuracy in one- and two-dimensional
problems. But it also tends to undercover and no accurate first order
version of this EEL is available. These motivated us to consider an EEL
defined on the full parameter space.

\section{Extended empirical likelihood on the full parameter space}

Consider a bijective mapping from the OEL domain to the parameter
space, $h_n\dvtx  \Theta_n \rightarrow\Theta=\mathbb{R}^d$. Under such a
mapping, the EEL $l^*(\theta)$ given by (\ref{1-5}) is well defined
throughout $\mathbb{R}^d$ and is thus free from the mismatch problem.
In this section, we first construct such a mapping using $h_n^*$ in
(\ref{2-80}). We call it the \emph{composite similarity mapping} and
denote it by $h^C_n\dvtx  \Theta_n\rightarrow\mathbb{R}^d$. We then study
the asymptotic properties of the EEL $l^*(\theta)$ based on $h^C_n$.

\subsection{The composite similarity mapping}

\begin{figure}

\includegraphics{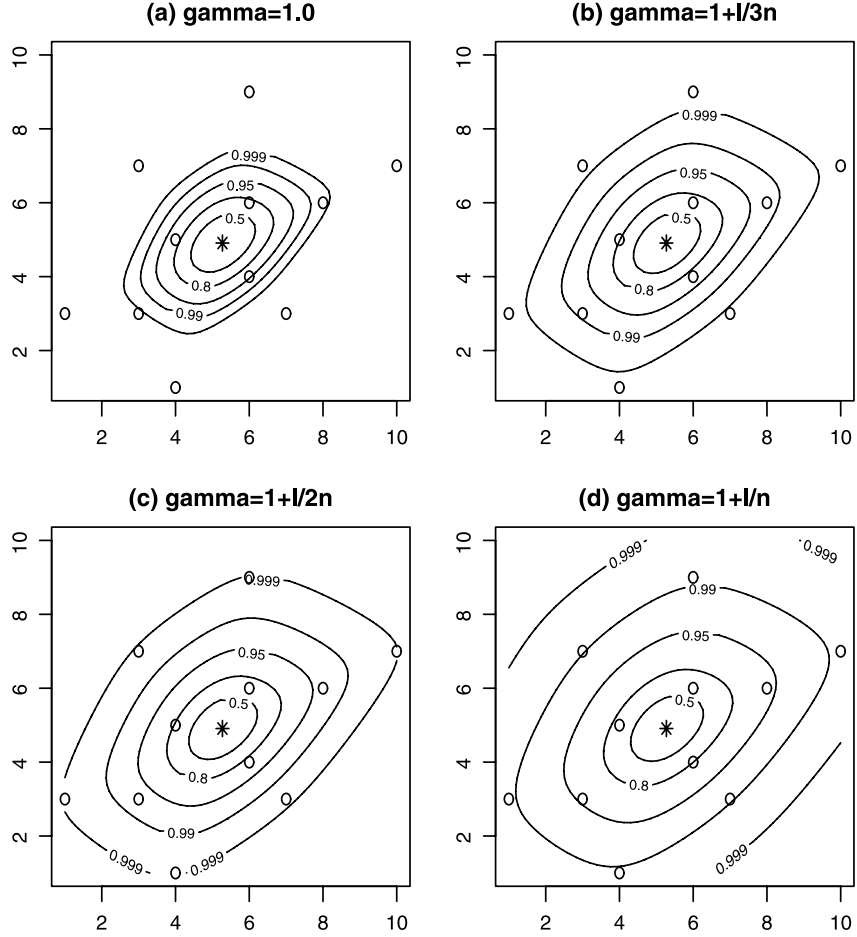}

\caption{\textup{(a)} Contours of the OEL $l(\theta)$ (for which the expansion
factor is 1.0); \textup{(b)} contours of the EEL $l^*(\theta)$ with expansion
factor $\gamma(11,l(\theta))=1+l(\theta)/3n$; \textup{(c)} contours of
the EEL
$l^*(\theta)$ with $\gamma(11,l(\theta))=1+l(\theta)/2n$; \textup{(d)} contours
of the EEL $l^*(\theta)$ with $\gamma(11,l(\theta))=1+l(\theta)/n$. All
four plots are based on the same random sample of 11 points shown in
small circles. The star in the middle is the sample mean. The expansion
factor increases as we go from plot \textup{(b)} to plot \textup{(d)}, and
correspondingly the EEL contours also become bigger in scale from plot
\textup{(b)} to \textup{(d)}. But the centre and shapes of contours are
the same in all
plots.}\label{fig1}
\end{figure}

The simple similarity mapping $h_n^*$ in (\ref{2-80}) maps OEL domain
$\Theta_n$ onto a similar but bounded region in $\mathbb{R}^d$.
If we think of $\Theta_n$ as a region consisting of distinct and nested
contours of the OEL, then $h_n^*$ expands all contours with the same
constant expansion factor $\gamma_n$. In order to map $\Theta_n$ onto
the full $\mathbb{R}^d$, we need to expand contours on the outside more
and more so that the images of the contours will fill up the entire
$\mathbb{R}^d$. To achieve this, consider level-$\tau$ contour of the
OEL $l(\theta)$,
%
\begin{equation}\label{3-2}
c(\tau)=\bigl\{\theta\dvtx  \theta\in\Theta_n\mbox{ and }
l(\theta)=\tau\bigr\},
\end{equation}
where $\tau\geq0$. The contours form a partition of the OEL domain,
%
\begin{equation}\label{3-3}
\Theta_n = \bigcup_{\tau\in[0,+\infty)} c(\tau).
\end{equation}
In light of (\ref{3-3}), the centre of $\Theta_n$ is $c(0)=\{\tilde
{\theta}\}$ and the outwardness of a $c(\tau)$ with respect to the
centre is indexed by $\tau$; the larger the $\tau$ value, the more
outward $c(\tau)$ is. If we allow the expansion factor $\gamma_n$ to be
a continuous monotone increasing function of $\tau$ and allow $\gamma
_n$ to go to infinity when $\tau$ goes to infinity, then (such a
variation of) $h_n^*$ will map $\Theta_n$ onto $\mathbb{R}^d$.
Hence, we define the composite similarity mapping $h^C_n\dvtx  \Theta
_n\rightarrow\mathbb{R}^d$ as follows:
%
\begin{equation}\label{3-15}
h_n^C(\theta)=\tilde{\theta}+\gamma\bigl(n,l(\theta)
\bigr) (\theta-\tilde{\theta})\qquad\mbox{for $\theta\in
\Theta_n$},
\end{equation}
where $\gamma(n,l(\theta))$ is given by
%
\begin{equation}\label{3-16}
\gamma\bigl(n,l(\theta)\bigr)=1+\frac{l(\theta)}{2n}.
\end{equation}
Function $\gamma(n,l(\theta))$ is the new expansion factor which
depends continuously on $\theta$ through the value of $l(\theta)$ or
$\tau=l(\theta)$. For convenience, we will emphasis the dependence of
$\gamma(n,l(\theta))$ on $l(\theta)$ instead of $\theta$ or $\tau$.
This new expansion factor has the two desired properties discussed above:
%
\begin{eqnarray}
\label{3-e1}
& & \mbox{for a fixed $n$, if $l(\theta_1)<l(\theta
_2)$, then $\gamma\bigl(n,l(\theta_2)\bigr)<\gamma
\bigl(n,l(\theta_2)\bigr)$; and}
\\
\label{3-e2}
& & \mbox{for a fixed $n$, $\gamma\bigl(n,l(\theta)\bigr)\rightarrow
+\infty$ as
$l(\theta) \rightarrow+\infty$.}
\end{eqnarray}
The inclusion of the denominator $2n$ in (\ref{3-16}) ensures that the
expansion factor converges to 1, reflecting the fact that there is no
need for domain expansion for large sample sizes.
Also, the constant $2$ in the denominator provides extra adjustment to
the speed of expansion and may be replaced with other positive
constants (see Figure \ref{fig1}). We choose to use 2 here as the corresponding
$\gamma(n,l(\theta))$ in (\ref{3-16}) is asymptotically equivalent to a
likelihood based expansion factor $\gamma(n, L(\theta))=L(\theta
)^{-1/n}$ which we had first considered and was found to give accurate
numerical results. The definition of $\gamma(n, l(\theta))$ in (\ref
{3-16}) uses $l(\theta)$ instead of $L(\theta)$ because of convenience
for theoretical investigations.
A more general form of $\gamma(n,l(\theta))$ will be considered later.

Theorem \ref{thm1} below summarizes the key properties of the composite
similarity mapping $h_n^C$. Its proof and that of subsequent theorems
and lemmas may all be found in the \hyperref[app]{Appendix}.

\begin{theorem} \label{thm1}
Under assumption ($A_1$), the composite similarity mapping $h^C_n\dvtx
\Theta_n\rightarrow\mathbb{R}^d$ defined by (\ref{3-15}) and (\ref
{3-16}) satisfies:

\begin{longlist}
\item
it has a unique fixed point at the mean
$\tilde{\theta}$;

\item it is a similarity mapping for each individual
$c(\tau)$; and

\item it is a bijective mapping from $\Theta_n$ to
$\mathbb{R}^d$.
\end{longlist}
\end{theorem}

Because of (ii) above, we call $h_n^C$ the composite similarity
mapping as it may be viewed as a continuous sequence of simple
similarity mappings from $\mathbb{R}^d$ to $\mathbb{R}^d$ indexed by
$\tau=l(\theta)\in[0, +\infty)$. The ``$\tau$th'' mapping from this
sequence has expansion factor $\gamma(n,l(\theta))=\gamma(n,\tau)$. It
is just the simple similarity mapping $h_n^*$ in (\ref{2-80}) with
$\gamma_n=\gamma(n,\tau)$, and is used exclusively to map the ``$\tau
$th'' OEL contour $c(\tau)$. The latter has been implicitly built into
$h_n^C$ since for all $\theta\in c(\tau)$, $l(\theta)=\tau$ which
implies the corresponding expansion factor of $h_n^C$ is the constant
$\gamma(n,\tau)$ that defines the ``$\tau$th'' mapping.

It should be noted that $h_n^C$ is not a similarity mapping from
$\mathbb{R}^d$ to $\mathbb{R}^d$ itself due to the dependence of the
expansion factor $\gamma(n,l(\theta))$ on $\theta$ and its domain
$\Theta_n$ which is only a bounded subset of $\mathbb{R}^d$.

\subsection{The extended empirical likelihood under the composite
similarity mapping}

By Theorem \ref{thm1}, $h_n^C\dvtx  \Theta_n\rightarrow\mathbb{R}^d$ is
bijective. Hence, it has an inverse which we denote by $h_n^{-C}\dvtx
\mathbb{R}^d \rightarrow\Theta_n$. The EEL $l^*(\theta)$ under
$h_n^C$ is
\[
l^*(\theta)=l\bigl(h_n^{-C}(\theta)\bigr) \qquad
\mbox{for $\theta\in\mathbb{R}^d$},
\]
which is defined throughout $\mathbb{R}^d$.
The contours of $l^*(\theta)$ are larger in scale but have the same
centre and identical shape as that of OEL $l(\theta)$.
Figure \ref{fig1} compares their contours with a sample of 11 two-dimensional
observations. It shows that geometrically, mapping $h_n^C$ is anchored
at the sample mean $\tilde{\theta}$ as it is the fixed point of $h_n^C$
that is not moved. From this anchoring point, the mapping pushes
out/expands each OEL contour $c(\tau)$ proportionally in all directions
at an expansion factor of $\gamma(n,\tau)$ to form an EEL contour. The
boundary points of $\Theta_n$ are all pushed out to the infinity.

In the following, we will use $\theta'$ to denote the image of a
$\theta\in\mathbb{R}^d$ under the inverse transformation $h_n^{-C}$,
that is, $h_n^{-C}(\theta)=\theta' \in\Theta_n$.
Of particular interest is the image of the unknown true mean $\theta_0$,
%
\begin{equation}\label{image}
h_n^{-C}(\theta_0)=\theta_0'.
\end{equation}
Because the inverse transformation $h_n^{-C}$ does not have an analytic
expression, that for $\theta'_0$ is also not available. Nevertheless,
Lemma \ref{lem2} gives an asymptotic assessment on its distance to
$\theta_0$. The proof of Lemma \ref{lem2} will need Lemma \ref{lem1}
below which shows that inside $\Theta_n$, the OEL $l(\theta)$ is a
``monotone increasing'' function along each ray originating from the
mean $\tilde{\theta}$.

\begin{lemma} \label{lem1}
Under assumption ($A_1$), for a fixed point $\theta\in\Theta_n$ and
any value $\alpha\in[0,1]$ the OEL $l(\theta)$ satisfies
\[
l\bigl(\tilde{\theta}+\alpha(\theta-\tilde{\theta})\bigr)\leq l(\theta).
\]
\end{lemma}

Denote by $[\tilde{\theta}, \theta_0]$ the line segment connecting
$\tilde{\theta}$ and $\theta_0$. Lemma \ref{lem2} below shows that
$\theta_0'$ is on $[\tilde{\theta}, \theta_0]$ and is asymptotically
very close to $\theta_0$.

\begin{lemma} \label{lem2} Under assumption ($A_1$), point $\theta_0'$
defined by equation (\ref{image}) satisfies \textup{(i)} $\theta_0' \in[\tilde
{\theta}, \theta_0]$ and \textup{(ii)} $\theta'_0-\theta_0 = O_p(n^{-3/2})$.
\end{lemma}

Using $\theta_0'$, the EEL $l^*(\theta_0)$ can now be expressed as
%
\begin{equation}\label{3-50}
l^*(\theta_0)=l\bigl(h_n^{-C}(
\theta_0)\bigr)= l\bigl(\theta'_0\bigr)=l
\bigl(\theta_0+\bigl(\theta_0'-
\theta_0\bigr)\bigr).
\end{equation}
The following theorem gives the asymptotic distribution of $l^*(\theta_0)$:

\begin{theorem} \label{thm2}
Under assumption ($A_1$) and with the composite similarity mapping
$h_n^C$ defined by (\ref{3-15}) and (\ref{3-16}), the EEL $l^*(\theta)$
satisfies
%
\begin{equation}\label{3-58}
l^*(\theta_0) \stackrel{\mathcal{D}} {\longrightarrow}
\chi^2_d \qquad\mbox{as $n\rightarrow+\infty$}.
\end{equation}
\end{theorem}
The proof of Theorem \ref{thm2} is based on the observation that $\|
\theta_0'-\theta_0\|$ is asymptotically very small. This and (\ref
{3-50}) imply that $l^*(\theta_0)=l(\theta_0)+o_p(1)$. This proof
demonstrates an advantage of the EEL: it has the simple relationship
with the OLE shown in (\ref{3-50}) through which we make use of known
asymptotic properties of the OEL to study the EEL. Our derivation of a
second order EEL below further explores this advantage.

\subsection{Second order extended empirical likelihood}

It may be verified that Theorems \ref{thm1} and \ref{thm2} also hold
under any composite similarity mapping defined by (\ref{3-15}) and the
following general form of the expansion factor,
%
\begin{equation}\label{3-90}
\gamma\bigl(n,l(\theta)\bigr)=1+\frac{\kappa[l(\theta)]^{\delta
(n)}}{n^m},
\end{equation}
where $\kappa$ and $m$ are both positive constants and $\delta(n)$ is a
bounded function of $n$ satisfying $\varepsilon_n< \delta(n) \leq a$
for some constants $a>\varepsilon_n>0$.
The availability of a whole family of $\gamma(n,l(\theta))$ functions
for the construction of the EEL provides an opportunity to optimize our
choice of this function to achieve the second order accuracy. Theorem
\ref{thm3} below gives the optimal choice.

\begin{theorem} \label{thm3} Assume $(A_1)$ and $(A_2)$ hold and denote
by $l^*_s(\theta)$ the EEL under the composite similarity mapping (\ref
{3-15}) with expansion factor
%
\begin{equation}\label{3-80}
\gamma_s\bigl(n,l(\theta)\bigr)=1+\frac{b }{2n}\bigl[l(\theta)
\bigr]^{\delta(n)},
\end{equation}
where $\delta(n)=O(n^{-1/2})$ and $b$ is the Bartlett correction
constant in (\ref{2-68}). Then
%
\begin{equation}\label{3-80-1}
l^*_s(\theta_0)=l(\theta_0)
\bigl[1-bn^{-1}+O_p\bigl(n^{-3/2}\bigr)\bigr]
\end{equation}
and
%
\begin{equation}\label{3-81}
P\bigl(l_s^*(\theta_0)\leq c\bigr)=P\bigl(
\chi^2_d \leq c\bigr)+O\bigl(n^{-2}\bigr).
\end{equation}
\end{theorem}

In our subsequent discussions, we will refer to an
$l^*_s(\theta)$ defined by the expansion factor $\gamma_s(n,l(\theta))$
in (\ref{3-80}) as a \emph{second order EEL} on the full parameter
space. The EEL $l^*(\theta)$ defined by $\gamma(n, l(\theta))$ in (\ref
{3-16}) will henceforth be referred to as the \emph{first order EEL} on
the full parameter space.

The utility of the $\delta(n)$ in $\gamma_s(n,l(\theta))$ is to provide
an extra adjustment for the speed of the domain expansion which ensures
that $l^*_s(\theta)$ will behave asymptotically like the BEL and hence
will have the second order accuracy of the BEL. For convenience, we set
$\delta(n)=n^{-1/2}$. The resulting second order EEL turns out to be
competitive in accuracy to the BEL and the second order AEL. For small
sample and/or high dimension situations, confidence regions based on
this second order EEL can have undercoverage problems like those based
on the OEL and BEL. Fine-tuning of $\delta(n)$ for such situations is
needed and methods of fine-tuning are discussed in \citet{Wu13}.

Finally, we noted after Theorem \ref{thm2} that the first order EEL
$l^*(\theta_0)$ can be expressed in terms of the OLE $l(\theta_0)$ as
$l^*(\theta_0)=l(\theta_0)+o_p(1)$. The $o_p(1)$ term can be improved
and in fact we have $l^*(\theta_0)=l(\theta_0)+O_p(n^{-1})$. See the
proof of Theorem \ref{thm2} in the \hyperref[app]{Appendix}. An even stronger
connection between $l^*(\theta_0)$ and $l(\theta_0)$ is given by
Corollary \ref{cor1} below.

\begin{cor} \label{cor1} Under assumptions ($A_1$) and ($A_2$), the
first order EEL $l^*(\theta)$ satisfies
%
\begin{equation}\label{3-81-1}
l^*(\theta_0)=l(\theta_0)\bigl[1-l(
\theta_0)n^{-1}+O_p\bigl(n^{-3/2}
\bigr)\bigr].
\end{equation}
\end{cor}

The proof of Corollary \ref{cor1} follows from that for Theorem \ref
{thm3}. This result provides a partial explanation for the remarkable
numerical accuracy of confidence regions based on the first order EEL
$l^*(\theta)$.

\begin{table}
\tabcolsep=0pt
\caption{Simulated coverage probabilities for
one-dimensional examples}\label{tab1}
{\fontsize{8.6pt}{11pt}\selectfont{
\begin{tabular*}{\tablewidth}{@{\extracolsep{\fill}}l c c c c c d{1.4} c c c c@{}}
\hline
& $\bolds{n}$ & \textbf{Level} & \textbf{OEL} & \textbf{\textit{EEL}$\bolds{_1}$}
& \multicolumn{1}{c}{\textbf{BEL}} & \multicolumn{1}{c}{\textbf{AEL}}
& \multicolumn{1}{c}{\textbf{EEL}$\bolds{_2}$} & \multicolumn{1}{c}{\textbf{BEL}$\bolds{^*}$} &
\multicolumn{1}{c}{\textbf{AEL}$\bolds{^*}$} & \multicolumn{1}{c@{}}{\textbf{EEL}$\bolds{^*_2}$} \\ \hline
$N(0,1)$
&10 &0.90 &0.8506 &\textbf{0.8914} &0.8753 &0.8788 &0.8813 &0.8767 &0.8867
&0.8824 \\
& &0.95 &0.9039 &\textbf{0.9452} &0.9246 &0.9294 &0.9317 &0.9242 &0.9352
&0.9324 \\
& &0.99 &0.9580 &\textbf{0.9867} &0.9677 &0.9753 &0.9738 &0.9656 &0.9771
&0.9734 \\
&30 &0.90 &0.8920 &\textbf{0.9071} &0.9007 &0.9008 &0.9022 &0.9017 &0.9019
&0.9030 \\
& &0.95 &0.9398 &\textbf{0.9548} &0.9461 &0.9461 &0.9476 &0.9466 &0.9468
&0.9474 \\
& &0.99 &0.9866 &\textbf{0.9925} &0.9882 &0.9883 &0.9885 &0.9883 &0.9884
&0.9884 \\
&50 &0.90 &0.8941 &\textbf{0.9024} &0.8995 &0.8996 &0.9003 &0.8992 &0.8993
&0.9000 \\
& &0.95 &0.9447 &\textbf{0.9541} &0.9481 &0.9479 &0.9486 &0.9483 &0.9484
&0.9490 \\
& &0.99 &0.9880 &\textbf{0.9920} &0.9892 &0.9892 &0.9892 &0.9894 &0.9894
&0.9895 \\
[4pt]
$t_5$
&10 &0.90 &0.8277 &\textbf{0.8765} &0.9226 &0.9979 &0.9209 &0.8520 &0.8873
&0.8782 \\
& &0.95 &0.8882 &\textbf{0.9394} &0.9599 &1.000 &0.9651 &0.9036 &0.9367
&0.9307 \\
& &0.99 &0.9556 &\textbf{0.9851} &0.9820 &1.000 &0.9887 &0.9499 &0.9798
&0.9751 \\
&30 &0.90 &0.8690 &\textbf{0.8852} &0.8999 &0.9028 &0.9017 &0.8852 &0.8885
&0.8882 \\
& &0.95 &0.9265 &\textbf{0.9436} &0.9476 &0.9509 &0.9502 &0.9385 &0.9428
&0.9420 \\
& &0.99 &0.9797 &\textbf{0.9888} &0.9875 &0.9901 &0.9886 &0.9831 &0.9863
&0.9861 \\
&50 &0.90 &0.8862 &\textbf{0.8967} &0.9040 &0.9048 &0.9052 &0.8977 &0.8983
&0.8987 \\
& &0.95 &0.9410 &\textbf{0.9491} &0.9515 &0.9518 &0.9518 &0.9465 &0.9471
&0.9474 \\
& &0.99 &0.9861 &\textbf{0.9918} &0.9907 &0.9913 &0.9913 &0.9881 &0.9882
&0.9886 \\[4pt]
$\chi_1^2$
&10 &0.90 &0.7764 &\textbf{0.8174} &0.8726 &1.000 &0.8634 &0.6792 &0.8456
&0.8291 \\
& &0.95 &0.8314 &\textbf{0.8781} &0.9068 &1.000 &0.9030 &0.7239 &0.8918
&0.8779 \\
& &0.99 &0.8973 &\textbf{0.9378} &0.9417 &1.000 &0.9461 &0.7677 &0.9391
&0.9253 \\
&30 &0.90 &0.8594 &\textbf{0.8759} &0.8887 &0.8901 &0.8890 &0.8658 &0.8847
&0.8829 \\
& &0.95 &0.9115 &\textbf{0.9249} &0.9319 &0.9343 &0.9330 &0.9105 &0.9278
&0.9272 \\
& &0.99 &0.9659 &\textbf{0.9764} &0.9759 &0.9786 &0.9769 &0.9565 &0.9735
&0.9733 \\
&50 &0.90 &0.8722 &\textbf{0.8833} &0.8936 &0.8941 &0.8943 &0.8887 &0.8912
&0.8909 \\
& &0.95 &0.9318 &\textbf{0.9411} &0.9441 &0.9458 &0.9459 &0.9388 &0.9419
&0.9415 \\
& &0.99 &0.9779 &\textbf{0.9847} &0.9837 &0.9845 &0.9845 &0.9804 &0.9831
&0.9830 \\[4pt]
$0.3N(0,1)$
&10 &0.90 &0.8470 &\textbf{0.8908} &0.8551 &0.8556 &0.8569 &0.8761 &0.8826
&0.8821 \\
$\quad{}+0.7N(2,1)$ & &0.95 &0.9036 &\textbf{0.9433} &0.9094 &0.9097 &0.9127 &0.9215 &0.9299
&0.9285 \\
& &0.99 &0.9564 &\textbf{0.9867} &0.9592 &0.9601 &0.9631 &0.9657 &0.9760
&0.9741 \\
&30 &0.90 &0.8930 &\textbf{0.9054} &0.8956 &0.8956 &0.8960 &0.9016 &0.9013
&0.9017 \\
& &0.95 &0.9438 &\textbf{0.9582} &0.9455 &0.9455 &0.9460 &0.9501 &0.9501
&0.9507 \\
& &0.99 &0.9873 &\textbf{0.9943} &0.9883 &0.9883 &0.9884 &0.9901 &0.9901
&0.9909 \\
&50 &0.90 &0.8965 &\textbf{0.9048} &0.8989 &0.8988 &0.8990 &0.9014 &0.9014
&0.9016 \\
& &0.95 &0.9465 &\textbf{0.9556} &0.9475 &0.9476 &0.9477 &0.9494 &0.9496
&0.9499 \\
& &0.99 &0.9876 &\textbf{0.9911} &0.9879 &0.9877 &0.9879 &0.9883 &0.9883
&0.9886 \\
\hline
\end{tabular*}}}
\end{table}

\section{Numerical examples and comparisons}\label{sec4}

We now present a simulation study comparing the EEL with the OEL, BEL
and AEL.
Throughout this section, we use $l^*_1(\theta)$ or EEL$_1$ to denote
the first order EEL with expansion factor (\ref{3-16}), and use
$l^*_2(\theta)$ or EEL$_2$ to denote the second order EEL given by
expansion factor (\ref{3-80}) where $\delta(n)=n^{-1/2}$.

\subsection{Low-dimensional examples}

Tables \ref{tab1} and \ref{tab2} contain simulated coverage probabilities of confidence
regions for the mean based on first order methods OEL, EEL$_1$ and
second order methods BEL, AEL, EEL$_2$, BEL$^*$, AEL$^*$ and EEL$^*_2$. Here,
BEL, AEL and EEL$_2$ are based on the theoretical Bartlett correction
constant $b$, and BEL$^*$, AEL$^*$ and EEL$^*_2$ are based on $\tilde{b}_n$
which is a bias corrected estimate for $b$ given by \citet{LiuChe10}.

%

\begin{table}
\tabcolsep=0pt
\caption{Simulated coverage probabilities for
two-dimensional examples}\label{tab2}
{\fontsize{8.6pt}{11pt}\selectfont{
\begin{tabular*}{\tablewidth}{@{\extracolsep{\fill}}l c c c c c d{1.4} c c c c@{}}
\hline
& $\bolds{n}$ & \textbf{Level} & \textbf{OEL} & \textbf{\textit{EEL}$\bolds{_1}$}
& \multicolumn{1}{c}{\textbf{BEL}} & \multicolumn{1}{c}{\textbf{AEL}}
& \multicolumn{1}{c}{\textbf{EEL}$\bolds{_2}$} & \multicolumn{1}{c}{\textbf{BEL}$\bolds{^*}$} &
\multicolumn{1}{c}{\textbf{AEL}$\bolds{^*}$} & \multicolumn{1}{c@{}}{\textbf{EEL}$\bolds{^*_2}$} \\
\hline
$\mathit{BV}1$
&10 &0.90 &0.7134 &\textbf{0.8118} &0.7965 &0.9777 &0.8212 &0.7561 &0.8350
&0.7989 \\
& &0.95 &0.7717 &\textbf{0.8758} &0.8407 &1.000 &0.8718 &0.8000 &0.8941
&0.8478 \\
& &0.99 &0.8484 &\textbf{0.9422} &0.8945 &1.000 &0.9268 &0.8570 &0.9680
&0.9083 \\
&30 &0.90 &0.8549 &\textbf{0.8888} &0.8824 &0.8856 &0.8872 &0.8786 &0.8813
&0.8822 \\
& &0.95 &0.9120 &\textbf{0.9426} &0.9313 &0.9348 &0.9361 &0.9296 &0.9336
&0.9337 \\
& &0.99 &0.9689 &\textbf{0.9868} &0.9772 &0.9798 &0.9796 &0.9757 &0.9783
&0.9779 \\
&50 &0.90 &0.8699 &\textbf{0.8917} &0.8869 &0.8874 &0.8894 &0.8859 &0.8869
&0.8886 \\
& &0.95 &0.9259 &\textbf{0.9428} &0.9354 &0.9361 &0.9374 &0.9344 &0.9351
&0.9363 \\
& &0.99 &0.9806 &\textbf{0.9908} &0.9846 &0.9852 &0.9856 &0.9839 &0.9848
&0.9851 \\[4pt]
$\mathit{BV}2$
&10 &0.90 &0.7513 &\textbf{0.8521} &0.8035 &0.8573 &0.8282 &0.7942 &0.8451
&0.8229 \\
& &0.95 &0.8061 &\textbf{0.9095} &0.8627 &0.9397 &0.8861 &0.8499 &0.9103
&0.8833 \\
& &0.99 &0.8879 &\textbf{0.9693} &0.9202 &1.000 &0.9430 &0.9116 &0.9721
&0.9405 \\
&30 &0.90 &0.8714 &\textbf{0.9019} &0.8864 &0.8872 &0.8897 &0.8864 &0.8881
&0.8906 \\
& &0.95 &0.9256 &\textbf{0.9549} &0.9406 &0.9413 &0.9428 &0.9403 &0.9409
&0.9432 \\
& &0.99 &0.9789 &\textbf{0.9907} &0.9823 &0.9826 &0.9838 &0.9820 &0.9829
&0.9836 \\
&50 &0.90 &0.8826 &\textbf{0.9037} &0.8935 &0.8937 &0.8954 &0.8938 &0.8939
&0.8952 \\
& &0.95 &0.9348 &\textbf{0.9528} &0.9423 &0.9426 &0.9438 &0.9423 &0.9425
&0.9435 \\
& &0.99 &0.9839 &\textbf{0.9914} &0.9862 &0.9864 &0.9871 &0.9861 &0.9864
&0.9864 \\[4pt]
$\mathit{BV}3$
&10 &0.90 &0.7001 &\textbf{0.7979} &0.7979 &1.000 &0.8162 &0.7333 &0.8363
&0.7911 \\
& &0.95 &0.7608 &\textbf{0.8581} &0.8374 &1.000 &0.8624 &0.7765 &0.8922
&0.8375 \\
& &0.99 &0.8331 &\textbf{0.9263} &0.8817 &1.000 &0.9151 &0.8286 &0.9639
&0.8942 \\
&30 &0.90 &0.8429 &\textbf{0.8775} &0.8749 &0.8789 &0.8788 &0.8719 &0.8780
&0.8764 \\
& &0.95 &0.9015 &\textbf{0.9363} &0.9266 &0.9326 &0.9319 &0.9221 &0.9282
&0.9280 \\
& &0.99 &0.9648 &\textbf{0.9817} &0.9740 &0.9776 &0.9760 &0.9709 &0.9750
&0.9740 \\
&50 &0.90 &0.8619 &\textbf{0.8836} &0.8807 &0.8820 &0.8836 &0.8787 &0.8801
&0.8808 \\
& &0.95 &0.9212 &\textbf{0.9403} &0.9351 &0.9364 &0.9379 &0.9325 &0.9346
&0.9347 \\
& &0.99 &0.9758 &\textbf{0.9848} &0.9810 &0.9816 &0.9817 &0.9802 &0.9806
&0.9810 \\[4pt]
$\mathit{BV}4$
&10 &0.90 &0.6408 &\textbf{0.7371} &0.7882 &1.000 &0.7940 &0.6240 &0.8382
&0.7596 \\
& &0.95 &0.7030 &\textbf{0.8027} &0.8212 &1.000 &0.8377 &0.6637 &0.8896
&0.8051 \\
& &0.99 &0.7788 &\textbf{0.8808} &0.8580 &1.000 &0.8914 &0.7129 &0.9576
&0.8602 \\
&30 &0.90 &0.8229 &\textbf{0.8595} &0.8709 &0.8820 &0.8760 &0.8598 &0.8738
&0.8681 \\
& &0.95 &0.8857 &\textbf{0.9191} &0.9215 &0.9329 &0.9255 &0.9079 &0.9212
&0.9170 \\
& &0.99 &0.9520 &\textbf{0.9734} &0.9689 &0.9819 &0.9717 &0.9591 &0.9696
&0.9667 \\
&50 &0.90 &0.8494 &\textbf{0.8707} &0.8758 &0.8783 &0.8783 &0.8716 &0.8755
&0.8740 \\
& &0.95 &0.9060 &\textbf{0.9251} &0.9256 &0.9287 &0.9282 &0.9221 &0.9259
&0.9251 \\
& &0.99 &0.9675 &\textbf{0.9807} &0.9781 &0.9797 &0.9793 &0.9753 &0.9778
&0.9765 \\
\hline
\end{tabular*}}}
\end{table}

Table \ref{tab1} gives four one-dimensional (1-$d$) examples. Table \ref{tab2} contains
four bivariate ($\mathit{BV}$ or 2-$d$) examples; the first three were taken
from \citet{LiuChe10}, the fourth is a ``2-$d$ chi-square'', and
here are the details:

($\mathit{BV}_1$): $X_1|D\sim N(0,D^2)$ and
$X_2|D\sim$ Gamma$(D^{-1},1)$.

($\mathit{BV}_2$): $X_1|D\sim$ Poisson$(D)$ and
$X_2|D\sim$ Poisson$(D^{-1})$.

($\mathit{BV}_3$): $X_1|D\sim$ Gamma$(D,1)$ and
$X_2|D\sim$ Gamma$(D^{-1},1)$.

($\mathit{BV}_4$): $X_1$ and $X_2$ are
independent copies of a $\chi^2_1$ random variable.

The $D$ in $\mathit{BV}_1$, $\mathit{BV}_2$ and $\mathit{BV}_3$ is a uniform random variable on
$[1,2]$ which is used to induce dependence between $X_1$ and $X_2$. We
included $n=10, 30, 50$ representing, respectively, small, medium and
large sample sizes. Each entry in the tables is based on 10,000 random
samples of size $n$, shown in column 2, from the distribution in column
1. Here are our observations:

\begin{longlist}[(2)]
\item[(1)] \textit{BEL}, \textit{AEL and EEL$_2$}: For $n=30$ and $50$,
these three theoretical second order methods are extremely close in
terms of coverage accuracy. This is to be expected as their coverage
errors are all $O(n^{-2})$ which is very small for medium or large
sample sizes.

For $n=10$, the AEL statistic suffers from a boundedness problem
[\citet{EmeOwe09}] which may lead to trivial 100\% confidence
regions or inflated coverage probabilities. This explains the 1.000's
in various places in the AEL column and renders the AEL unsuitable for
such small sample sizes. Between BEL and EEL$_2$, the latter is more
accurate, especially for the 2-$d$ examples.

Overall, EEL$_2$ is the most accurate theoretical second order method.

\item[(2)] \textit{BEL$^*$}, \textit{AEL$^*$ and EEL$^*_2$}: For $n=30$ and $50$,
the AEL$^*$ and EEL$^*_2$ are slightly more accurate than the BEL$^*$,
especially in 2-$d$ examples.

For $n=10$, the AEL$^*$ has higher coverage probabilities but these are
inflated by and unreliable due to the boundedness problem. Also,
EEL$^*_2$ is more accurate than BEL$^*$. For the ``2-$d$ chi-square'' in
example $\mathit{BV}_4$, EEL$^*_2$ is at least 12\% more accurate than the BEL$^*$.

Overall, EEL$^*_2$ is the most reliable and accurate among the three.

\item[(3)] \textit{OEL and EEL$_1$}: These first order methods are
simpler than the second order methods as they do not require
computation of the theoretical or estimated Bartlett correction factor.
The EEL$_1$ is consistently and substantially more accurate than the
OEL. In particular, for 2-$d$ examples with $n=10$, the EEL$_1$ is more
accurate by about 10\%.

\begin{table}
\tabcolsep=0pt
\caption{Simulated coverage probabilities for the mean
of $d$-dimensional multivariate normal distributions}\label{tab3}
{\fontsize{8.6pt}{11pt}\selectfont{
\begin{tabular*}{\tablewidth}{@{\extracolsep{\fill}}l c c c c c d{1.4} c c d{1.4} c@{}}
\hline
& $\bolds{n}$ & \textbf{Level} & \textbf{OEL} & \textbf{\textit{EEL}$\bolds{_1}$}
& \multicolumn{1}{c}{\textbf{BEL}} & \multicolumn{1}{c}{\textbf{AEL}}
& \multicolumn{1}{c}{\textbf{EEL}$\bolds{_2}$} & \multicolumn{1}{c}{\textbf{BEL}$\bolds{^*}$} &
\multicolumn{1}{c}{\textbf{AEL}$\bolds{^*}$} & \multicolumn{1}{c@{}}{\textbf{EEL}$\bolds{^*_2}$} \\
\hline
$d=5$
&10 &0.90 &0.3007 &\textbf{0.5897} &0.3839 &1.000 &0.5306 &0.3691 &1.000
&0.5005 \\
& &0.95 &0.3368 &\textbf{0.6794} &0.4135 &1.000 &0.5842 &0.4028 &1.000
&0.5500 \\
& &0.99 &0.3946 &\textbf{0.7984} &0.4498 &1.000 &0.6642 &0.4422 &1.000
&0.6287 \\
&30 &0.90 &0.7790 &\textbf{0.8862} &0.8258 &0.8554 &0.8455 &0.8273 &0.8585
&0.8468 \\
& &0.95 &0.8497 &\textbf{0.9436} &0.8880 &0.9208 &0.9047 &0.8884 &0.9256
&0.9052 \\
& &0.99 &0.9337 &\textbf{0.9881} &0.9532 &0.9889 &0.9629 &0.9539 &0.9899
&0.9634 \\
&50 &0.90 &0.8476 &\textbf{0.9036} &0.8752 &0.8803 &0.8820 &0.8757 &0.8808
&0.8825 \\
& &0.95 &0.9089 &\textbf{0.9522} &0.9297 &0.9341 &0.9349 &0.9297 &0.9349
&0.9354 \\
& &0.99 &0.9728 &\textbf{0.9913} &0.9804 &0.9833 &0.9830 &0.9804 &0.9839
&0.9831 \\[4pt]
$d=10$
&20 &0.90 &0.1889 &\textbf{0.5367} &0.2845 &1.000 &0.4297 &0.2823 &1.000
&0.4235 \\
& &0.95 &0.2281 &\textbf{0.6260} &0.3209 &1.000 &0.4905 &0.3191 &1.000
&0.4824 \\
& &0.99 &0.2895 &\textbf{0.7708} &0.3747 &1.000 &0.5783 &0.3727 &1.000
&0.5717 \\
&30 &0.90 &0.4689 &\textbf{0.7752} &0.5944 &1.000 &0.6750 &0.5954 &1.000
&0.6752 \\
& &0.95 &0.5432 &\textbf{0.8594} &0.6627 &1.000 &0.7492 &0.6635 &1.000
&0.7480 \\
& &0.99 &0.6698 &\textbf{0.9442} &0.7670 &1.000 &0.8514 &0.7675 &1.000
&0.8527 \\
&50 &0.90 &0.7097 &\textbf{0.8806} &0.7921 &0.9531 &0.8189 &0.7933 &0.9582
&0.8198 \\
& &0.95 &0.7959 &\textbf{0.9393} &0.8577 &0.9968 &0.8827 &0.8588 &0.9974
&0.8838 \\
& &0.99 &0.9027 &\textbf{0.9864} &0.9392 &1.000 &0.9546 &0.9396 &1.000
&0.9549 \\[4pt]
$d=15$
&30 &0.90 &0.1224 &\textbf{0.4850} &0.2199 &1.000 &0.3581 &0.2196 &1.000
&0.3569 \\
& &0.95 &0.1513 &\textbf{0.5761} &0.2504 &1.000 &0.4130 &0.2502 &1.000
&0.4124 \\
& &0.99 &0.2155 &\textbf{0.7490} &0.3100 &1.000 &0.5077 &0.3099 &1.000
&0.5054 \\
&50 &0.90 &0.4769 &\textbf{0.7983} &0.6177 &1.000 &0.6883 &0.6191 &1.000
&0.6894 \\
& &0.95 &0.5665 &\textbf{0.8776} &0.6971 &1.000 &0.7630 &0.6985 &1.000
&0.7646 \\
& &0.99 &0.7065 &\textbf{0.9600} &0.8097 &1.000 &0.8682 &0.8103 &1.000
&0.8686 \\
&100 &0.90 &0.7696 &\textbf{0.9031} &0.8325 &0.9309 &0.8472 &0.8328 &0.9341
&0.8484 \\
& &0.95 &0.8484 &\textbf{0.9514} &0.8985 &0.9852 &0.9086 &0.8989 &0.9865
&0.9096 \\
& &0.99 &0.9405 &\textbf{0.9900} &0.9639 &1.000 &0.9693 &0.9641 &1.000
&0.9692 \\
\hline
\end{tabular*}}}
\end{table}

\item[(4)] \textit{EEL$_1$ versus EEL$^*_2$}: These are the most
accurate practical first and second order methods, respectively.
Surprisingly, EEL$_1$ turns out to be slightly more accurate than
EEL$^*_2$. Only the (impractical) theoretical second order EEL$_2$ is
comparable to EEL$_1$ in accuracy. This intriguing observation may be
partially explained by Corollary \ref{cor1} where it was shown that
$l^*_1(\theta_0)=l(\theta_0)[1+l(\theta_0)n^{-1}+O(n^{-3/2})]$, which
resembles the Bartlett corrected OEL in (\ref{2-72}) with the constant
$b$ replaced by $l(\theta_0)$. However, this does not account for its
good accuracy for small sample sizes, which is due to the fact that
EEL$_1$ makes good use of the dimension information through the
composite similarity mapping. We will further elaborate on this in
Section~\ref{sec4.2}.

\item[(5)] \textit{EEL$_1$}: Overall, it is the most accurate
among the eight methods that we have compared. Importantly, it is not
just accurate in relative terms. It is sufficiently accurate in
absolute terms for practical applications in most 1-$d$ examples,
including cases of $n=10$. It is also quite accurate for 2-$d$ examples
when $n=30, 50$.
\end{longlist}

\subsection{High-dimensional examples}\label{sec4.2}

Table \ref{tab3} contains simulated coverage probabilities for the mean of three
high-dimensional multivariate normal distributions ($d=5,10,15$). Our
main interest here is to probe the small sample behaviour of all
methods in high-dimension situations. Because of this, we have included
only combinations of $n$ and $d$ where $n/d$, which we will refer to as
the \emph{effective sample size}, is very small ($2\leq n/d\leq10$).
The following are our observations based on Table \ref{tab3}:
\begin{longlist}[(2)]
\item[(1)] For these high-dimension examples, EEL$_1$ is the most
accurate, surpassing even the theoretical second order EEL$_2$. Whereas
the OEL uses dimension $d$ only once through the degrees of freedom in
the chi-square calibration, EEL$_1$ uses $d$ twice. The expansion
factor for EEL$_1$ is $1+l(\theta)/2n$ which implicitly depends on $d$;
the 100($1-\alpha$)\% EEL$_1$ confidence region is just the
100($1-\alpha$)\% OEL confidence region expanded by a factor of $1+\chi
^2_{d,1-\alpha}/2n$. Hence, EEL$_1$ uses $d$ through the chi-square
calibration of the OEL region {and} the expansion factor.

%
\begin{table}[b]
\caption{Average lengths of EL confidence intervals for
$N(0,1)$ mean} \label{tab4}
\begin{tabular*}{\tablewidth}{@{\extracolsep{\fill}}l c c c c c c@{}}
\hline
$\bolds{n}$ &\textbf{Level}& \textbf{OEL} & \textbf{EEL}$\bolds{_1}$ &
\textbf{BEL}$\bolds{^*}$ & \textbf{AEL}$\bolds{^*}$ & \textbf{EEL}$\bolds{^*_2}$ \\
\hline
10 &0.90 &0.965 &1.096 &1.044 &$N/A$ (0.026)
&1.077 \\
&0.95 &1.149 &1.370 &1.242 &$N/A$ (0.058)
&1.298 \\
&0.99 &1.499 &1.996 &1.615 &$N/A$ (0.172)
&1.731 \\[4pt]
30 &0.90 &0.589 &0.616 &0.606 &0.606 &0.608 \\
&0.95 &0.706 &0.752 &0.726 &0.727 &0.731 \\
&0.99 &0.940 &1.044 &0.967 &0.969 &0.976 \\[4pt]
50 &0.90 &0.460 &0.473 &0.467 &0.468 &0.468 \\
&0.95 &0.551 &0.572 &0.560 &0.560 &0.561 \\
&0.99 &0.732 &0.780 &0.744 &0.744 &0.747 \\ \hline
\end{tabular*}
\end{table}

For a fixed $\alpha\in(0,1)$, the chi-square quantile $\chi
^2_{d,1-\alpha}$ and consequently the EEL$_1$ expansion factor $1+\chi
^2_{d,1-\alpha}/2n$ are increasing functions of $d$. Hence at a fixed
$n$, EEL$_1$ automatically provides higher degrees of expansion for
higher dimensions where this is indeed needed.

\item[(2)] For multivariate normal means, Table \ref{tab3} shows that EEL$_1$ is
accurate when the effective sample size satisfies $n/d\geq5$. However,
when the underlying distribution is heavily skewed, the effective
sample size needed to achieve similar accuracy needs to be 15 or
larger. See Table \ref{tab2} for some 2-$d$ examples to this effect.

\item[(3)] The AEL and AEL$^*$ broke down in most cases with 100\%
coverage probabilities. This further illustrates the observation that
AEL methods may not be suitable when the effective sample size is
small. Among OEL, BEL, EEL$_2$ and EEL$^*_2$, the two EEL methods are
consistently more accurate but they are not sufficiently accurate for
practical applications except for the case of $(d,n)=(5,50)$.
\end{longlist}

\subsection{Confidence region size comparison} For the 1-$d$ examples
in Table~\ref{tab1}, we computed the average interval lengths of the five
practical methods OEL, EEL$_1$, BEL$^*$, AEL$^*$ and EEL$^*_2$. Table \ref{tab4} gives
the average length of 1000 intervals of each method and $n$ combination
for the $N(0,1)$ case. For $n=10$, the average for AEL$^*$ is not
available due to occurrences of unbounded intervals; the number beside
the $N/A$ is 
the proportion of times where this occurred.
Not surprisingly, intervals with higher coverage probabilities in
Table \ref{tab1} have larger average lengths. That of EEL$_1$ is the largest but
it is not excessive relative to averages of other methods. As such,
length is not a big disadvantage for EEL$_1$ as other methods have
lower coverage probabilities.

For $d>1$, sizes of the confidence regions are difficult to determine.
But the relative size of an EEL region to the corresponding OEL region
can be measured by the expansion factor. Table \ref{tab5} contains values of the
expansion factor for 95\% EEL$_1$ regions at some combinations of $n$
and $d$.
The expansion factor increases when $d$ goes up but decreases when $n$
goes up, responding to the need for more expansion in higher dimension
situations and the need for less expansion when the sample size is large.

%
\begin{table}
\caption{Values of expansion factor for 95\% EEL$_1$
confidence regions}\label{tab5}
\begin{tabular*}{\tablewidth}{@{\extracolsep{\fill}}l c c c c c c@{}}
\hline
$\bolds{n}$ &$\bolds{d=1}$ & $\bolds{d=2}$ & $\bolds{d=3}$ & $\bolds{d=5}$
& $\bolds{d=10}$ & $\bolds{d=15}$ \\
\hline
10 &1.192 &1.299 &1.390 &1.553 &$N/A$ & $N/A$ \\
15 &1.128 &1.199 &1.260 &1.369 &1.610 & $N/A$ \\
20 &1.096 &1.149 &1.195 &1.276 &1.457 &1.624 \\
30 &1.064 &1.099 &1.130 &1.184 &1.305 &1.416 \\
50 &1.038 &1.059 &1.078 &1.110 &1.183 &1.249\\
\hline
\end{tabular*}
\end{table}

Finally, we briefly comment on the computation concerning the EEL
$l^*_1(\theta)$. To compute $l^*_1(\theta)$ at a given $\theta\in
\mathbb{R}^d$ which is just OEL $l(\theta')$ where $\theta'$ satisfies
equation $h_n^{C}(\theta')=\theta$, we need to find the multivariate
root for function $f(\theta')=h_n^{C}(\theta')-\theta$. This is seen as
a nonlinear multivariate problem but it is easily reduced to a simpler
\emph{univariate} problem due to the fact that $\theta' \in[\tilde
{\theta}, \theta]$ (see Lemma \ref{thm2} and its proof). When using
$l^*_1(\theta)$ for hypothesis testing or when simulating the coverage
probabilities of the EEL confidence regions, we may use the fact that
$l^*_1(\theta)\leq l(\theta)$. Hence, we can compute $l(\theta)$ first
and if it is smaller than the critical value, then there is no need to
compute $l^*_1(\theta)$ because it must also be smaller than the
critical value.
Incorporating these observations, our R code for computing the EEL runs
quite fast.

%
\section{Concluding remarks}

The geometric motivation of the domain expansion method is simple:
since the OEL confidence region tends to be too small, an expansion of
the OEL confidence region should help to ease its undercoverage
problem. 
What needed to be determined then are the manner in which the expansion
should take place and the amount of expansion that would be
appropriate. The composite similarity mapping of the present paper is
an effective way to undertake the expansion as it solves the mismatch
problem and retains all important geometric characteristics of the OEL
contours. With the impressive numerical accuracy of the EEL$_1$, the
particular amount of expansion represented by its expansion factor (\ref
{3-16}) would be appropriate for general applications of the EEL method.

The EEL is readily constructed for parameters defined by general
estimating equations. For such parameters, we use the maximum empirical
likelihood estimator (MELE) $\tilde{\theta}$ to define the composite
similarity mapping in (\ref{3-15}). Under certain conditions on the
estimating function which also guarantee the $\sqrt{n}$-consistency of
the MELE, Lemma \ref{thm2} and all three theorems of this paper remain
valid. A detailed treatment of the EEL in this setting may be found in
a technical report by \citet{TsaWu13}. See also \citet{Tsa13} for an
EEL for estimating equations under the simple similarity transformation
(\ref{2-80}). For parameters outside of the standard estimating
equations framework, the EEL on full parameter space may also be
defined through a composite similarity mapping centred on the MELE, but
its asymptotic properties need to be investigated for each case separately.

To conclude, the simple first order EEL$_1$ is a practical and reliable
method that is remarkably accurate when the effective sample size is
not too small. It is also easy to use. Hence, we recommend it for real
applications of the empirical likelihood method.
An intriguing question that remains largely unanswered is why this
first order method has such good accuracy relative to the OEL and the
second order methods. Corollary \ref{cor1} and the first remark in
Section \ref{sec4.2} suggested, respectively, possible asymptotic and finite
sample reasons, but a more convincing theoretical explanation is needed.

\begin{appendix}\label{app}
\section*{Appendix} We give proofs for lemmas and theorems below.

\begin{pf*}{Proof of Theorem \ref{thm1}}
Assumption ($A_1$) and $n>d$ imply that, with probability 1, the convex
hull of the $X_i$ is nondegenerate. This implies the OEL $l(\theta)$
has an open domain $\Theta_n\subset\mathbb{R}^d$, a condition that is
required for the implementation of OEL domain expansion via composite
similarity mapping. Subsequent proofs all require this condition which,
hereafter, is assumed whenever $(A_1)$ is, and for brevity will not be
explicitly restated.

Part (i) is a simple consequence of the observation that $\gamma
(n,l(\theta))\geq1$. To show~(ii), let $n$ and $\tau$ be fixed, and
consider the level-$\tau$ OEL contour $c(\tau)$ defined by (\ref{3-2}).
For $\theta\in c(\tau)$, $l(\theta)=\tau$. Thus, the composite
similarity mapping $h_n^C$ simplifies to
$h_n^C(\theta)=\tilde{\theta}+\gamma_n(\theta-\tilde{\theta})$ for
$\theta\in c(\tau)$
where $\gamma_n=\gamma(n,\tau)$ is constant. This is the simple
similarity mapping in (\ref{2-80}).

To prove (iii), we need to show that $h_n^C\dvtx  \Theta_n \rightarrow
\mathbb{R}^d$ is both surjective and injective. We first show it is
surjective, that is, for any given $\theta' \in\mathbb{R}^d$, there
exists a $\theta'' \in\Theta_n$ such that $h_n^C(\theta'')=\theta'$.
Consider the ray originating from $ \tilde{\theta}$ and through~$\theta
'$. Introduce a univariate parametrization of this ray,
\[
\theta= \theta(\zeta)= \tilde{\theta}+\zeta\vec{\theta},
\]
where $\vec{\theta}$ is the unit vector $(\theta'-\tilde{\theta})/\|
\theta'-\tilde{\theta} \|$ in the direction of the ray and $\zeta\in
[0, \infty)$ is the distance between $\theta$ (a point on the ray) and
$\tilde{\theta}$. Define
\[
\zeta_b=\inf\bigl\{\zeta\dvtx  \zeta\in[0,+\infty) \mbox{
and } \theta(\zeta) \notin\Theta_n \bigr\}.
\]
Then, $\theta(\zeta)\in\Theta_n$ for all $\zeta\in[0, \zeta_b)$. But
$\theta(\zeta_b)\notin\Theta_n$ because $\Theta_n$ is open. It follows
that $\zeta_b>0$ as it represents the distance between $\tilde{\theta
}$, an interior point of the open~$\Theta_n$, and $\theta(\zeta_b)$
which is a boundary point of $\Theta_n$.

Now, consider the following univariate function defined on $[0,
\zeta_b)$:
\[
f(\zeta) = \gamma\bigl(n,l\bigl(\theta(\zeta)\bigr)\bigr)\zeta. %
\]
We have $f(0)=\gamma(n,l(\tilde{\theta}))\times0=\gamma(n,0)\times
0=0$. Also, by (\ref{3-e2}),
\[
\lim_{\zeta\rightarrow\zeta_b} f(\zeta)= \lim_{\zeta\rightarrow\zeta
_b}\gamma
\bigl(n,l\bigl(\theta(\zeta)\bigr)\bigr) \zeta= \zeta_b \lim
_{\zeta\rightarrow\zeta_b}\gamma\bigl(n,l\bigl(\theta(\zeta)\bigr)\bigr
)=+\infty.
\]
Hence, by the continuity of $f(\zeta)$, for $\zeta'=\| \theta'-\tilde
{\theta}\| \in[0,+\infty)$, there exists a $\zeta'' \in[0,\zeta_b)$
such that $f(\zeta'') = \zeta'$. Let $\theta''=\theta(\zeta'')$. Then
$\theta'' \in\Theta_n$ since $\zeta'' \in[0,\zeta_b)$. Also,
$h_n^C(\theta'') = \tilde{\theta} + \gamma(n, l(\theta''))(\theta
''-\tilde{\theta})= \theta'$.
Hence, $\theta''$ is the desired point in $\Theta_n$ that satisfies
$h_n^C(\theta'')=\theta'$ and $h_n^C$ is surjective.

To show that $h_n^C$ is also injective, first note that for a given OEL
contour $c(\tau)$, the mapping $h_n^{C}\dvtx  c(\tau) \rightarrow c^*(\tau)$
is injective\vspace*{2pt} because for $\theta\in c(\tau)$, $h_n^C$ is equivalent to
the similarity mapping in (\ref{2-80}) which is bijective from $\mathbb
{R}^d$ to $\mathbb{R}^d$. By the partition of the OEL domain $\Theta_n$
in (\ref{3-3}),
two different points $\theta_1$, $\theta_2$ from $\Theta_n$ are either
[$a$] on the same contour $c(\tau)$ where $\tau=l(\theta_1)=l(\theta
_2)$ or [$b$] on two separate contours $c(\tau_1)$ and $c(\tau_2)$,
respectively, where $\tau_1=l(\theta_1)\neq l(\theta_2)=\tau_2$. Under
[$a$], $h_n^C(\theta_1)\neq h_n^C(\theta_2)$ because $h_n^{C}\dvtx  c(\tau)
\rightarrow c^*(\tau)$ is injective. Under [$b$], $h_n^C(\theta_1)\neq
h_n^C(\theta_2)$ also holds as (\ref{3-e1}) implies $c^*(\tau_1) \cap
c^*(\tau_2)=\varnothing$.
\end{pf*}

\begin{pf*}{Proof of Lemma \ref{lem1}}
For a fixed $\theta\in\Theta_n$, $l(\theta)$ is a fixed quantity in
$[0,+\infty)$. Define an OEL confidence region for the mean using
$l(\theta)$ as follows:
%
\begin{equation}\label{3-20}
{\cal C}_\theta=\bigl\{\theta'\dvtx  \mbox{ $
\theta' \in\mathbb{R}^d$ and $l\bigl(\theta'
\bigr) \leq l(\theta)$} \bigr\}.
\end{equation}
Then, ${\cal C}_\theta$ is a convex set in $\mathbb{R}^d$. See \citet{Owe90} and \citet{HalLaS90}. Since $l(\tilde{\theta})=0$ and
$l(\theta) \geq0$, $\tilde{\theta}$ is in ${\cal C}_\theta$. Further,
by construction, $\theta$ itself is also in~${\cal C}_\theta$. It
follows from the convexity of
${\cal C}_\theta$ that for any $\alpha\in[0,1]$,
\[
\theta^*=(1-\alpha)\tilde{\theta} + \alpha\theta= \tilde{\theta} +
\alpha(
\theta-\tilde{\theta})
\]
must also be in ${\cal C}_\theta$. By (\ref{3-20}), $l(\theta^*) \leq
l(\theta)$. Thus, $ l(\tilde{\theta}+\alpha(\theta-\tilde{\theta
}))\leq l(\theta)$. 
\end{pf*}

\begin{pf*}{Proof of Lemma \ref{lem2}}
Since $\theta_0= h_n^C(\theta_0')=\tilde{\theta}+\gamma(n,l(\theta
_0'))(\theta_0'-\tilde{\theta})$,
%
\begin{equation}\label{3-25}
\theta_0-\tilde{\theta}=\gamma\bigl(n,l\bigl(\theta_0'
\bigr)\bigr) \bigl(\theta_0'-\tilde{\theta}\bigr).
\end{equation}
Noting that $\gamma(n,l(\theta))\geq1$, (\ref{3-25}) implies $\theta
'_0$ is on the ray originating from $ \tilde{\theta}$ and through
$\theta_0$ and $\|\theta_0-\tilde{\theta}\|\geq\|\theta_0'-\tilde{\theta
}\|$. Hence, $\theta_0'\in[\tilde{\theta}, \theta_0]$.

Without\vspace*{1pt} loss of generality, we assume that $\theta_0 \in\Theta_n$. See
\citet{Owe90}. By the convexity of $\Theta_n$, $[\tilde{\theta}, \theta
_0] \subset\Theta_n$. It follows from Lemma \ref{lem1} that
\[
0=l(\tilde{\theta}) \leq l\bigl(\theta_0'\bigr) \leq
l(\theta_0). %
\]
This and the fact that $l(\theta_0) = O_p(1)$ imply $l(\theta
_0')=O_p(1)$. Hence,
%
\begin{equation}\label{3-30}
\gamma\bigl(n,l\bigl(\theta_0'\bigr)\bigr)=1+
\frac{l(\theta_0')}{2n}=1+O_p\bigl(n^{-1}\bigr).
\end{equation}
Replacing $(\theta_0'-\tilde{\theta})$ in (\ref{3-25}) with $(\theta
_0'-\theta_0+\theta_0-\tilde{\theta})$, we obtain
%
\begin{equation}\label{3-34}
\bigl[1-\gamma\bigl(n,l\bigl(\theta_0'\bigr)\bigr)
\bigr](\theta_0-\tilde{\theta})=\gamma\bigl(n,l\bigl(\theta
_0'\bigr)\bigr) \bigl(\theta_0'-
\theta_0\bigr).
\end{equation}
By $\tilde{\theta}-\theta_0=O_p(n^{-1/2})$, (\ref{3-30}) and (\ref
{3-34}), we have $\theta_0'-\theta_0=O_p(n^{-3/2})$.
\end{pf*}

It may be verified using the same steps in the above proof that if the
expansion factor $\gamma(n,l(\theta))$ in (\ref{3-16}) is replaced with
a more general $\gamma(n,l(\theta))=1+O_p(n^{-m})$ such as that in (\ref
{3-90}) where $m>0$, then Lemma \ref{lem2}(i) still holds and (ii)
becomes $\theta_0'-\theta_0=O_p(n^{-m-1/2})$. In particular, under
expansion factor $\gamma_s(n,l(\theta))$ in (\ref{3-80}), we also have
$\theta_0'-\theta_0=O_p(n^{-3/2})$.

\begin{pf*}{Proof of Theorem \ref{thm2}}
Differentiating both sides of (\ref{2-58}), we obtain
$\partial l(\theta)/\partial\theta=-2n\lambda^T$. By (ii) in Lemma
\ref{lem2}, $\theta_0'-\theta_0=O_p(n^{-3/2})$. Applying Taylor's
expansion to $l^*(\theta_0)=l(\theta_0+(\theta_0'-\theta_0))$ in (\ref
{3-50}), we have
%
\begin{equation}\label{3-61}
l^*(\theta_0)=l(\theta_0)-2n\lambda(
\theta_0)^T\bigl(\theta_0'-
\theta_0\bigr) + o_p\bigl(\bigl\|\bigl(\theta_0'-
\theta_0\bigr)\bigr\|\bigr).
\end{equation}
By \citet{Owe90}, $\lambda(\theta_0)=O_p(n^{-1/2})$. This and (\ref
{3-61}) imply that
$ l^*(\theta_0)=l(\theta_0)+ O_p(n^{-1})$, which together with (\ref
{2-60}), imply Theorem \ref{thm2}.
\end{pf*}

For cases where $\theta_0'-\theta_0=o_p(n^{-1/2})$, we have $l^*(\theta
_0)=l(\theta_0)+ o_p(1)$ which also implies Theorem \ref{thm2}.
Since $\theta_0'-\theta_0=o_p(n^{-1/2})$ under expansion factor (\ref
{3-90}), Theorem \ref{thm2} also holds for EEL defined by expansion
factor (\ref{3-90}). 

\begin{pf*}{Proof of Theorem \ref{thm3}}
First, note that $\gamma_s(n,l(\theta))$ in (\ref{3-80}) satisfies
conditions (\ref{3-e1}) and (\ref{3-e2}). Thus it may be verified that
Theorem \ref{thm1}, Lemma~\ref{lem2} and Theorem \ref{thm2} hold under
the composite similarity mapping given by $\gamma_s(n,l(\theta))$.
In particular, the EEL corresponding to this composite similarity
mapping, $l^*_s(\theta_0)$, converges in distribution to a $\chi^2_d$
random variable.

Since $\delta(n)=O(n^{-1/2})$ and $l(\theta_0')=l^*_s(\theta
_0)=O_p(1)$, we have
%
\begin{equation}\label{3-83}
\bigl[l\bigl(\theta_0'\bigr)\bigr]^{\delta(n)}=1+O_p
\bigl(n^{-1/2}\bigr).
\end{equation}
With expansion factor $\gamma_s(n,l(\theta))$ in (\ref{3-80}), equation
(\ref{3-25}) becomes
\[
\theta_0-\tilde{\theta}=\gamma_s\bigl(n,l\bigl(
\theta_0'\bigr)\bigr) \bigl(\theta_0'-
\tilde{\theta}\bigr).
\]
This implies
%
\begin{eqnarray}\label{3-87-3}
\theta_0'-\theta_0&=&\frac{b [l(\theta_0')]^{\delta(n)}}{2n}
\bigl(\tilde{\theta}-\theta_0'\bigr)\nonumber\\[-8pt]\\[-8pt]
&=&\frac{b [l(\theta_0')]^{\delta(n)}}{2n}(
\tilde{\theta}-\theta_0)+\frac{b [l(\theta_0')]^{\delta(n)}}{2n}\bigl(
\theta_0-\theta_0'\bigr).\nonumber
\end{eqnarray}
It follows from (\ref{3-83}), (\ref{3-87-3}) and $\theta_0'-\theta
_0=O_p(n^{-3/2})$ that
%
\begin{equation}\label{3-87-2}
\theta_0'-\theta_0=\frac{b [l(\theta_0')]^{\delta(n)}}{2n}(
\tilde{\theta}-\theta_0)+O_p\bigl(n^{-5/2}
\bigr) =\frac{b}{2n}(\tilde{\theta}-\theta_0)+O_p
\bigl(n^{-2}\bigr).\hspace*{-25pt}
\end{equation}
By assumptions ($A_1$) and ($A_2$), the OEL $l(\theta_0)$ has expansion
%
\begin{equation}\label{3-85}
l(\theta_0) = n(\tilde{\theta} -\theta_0)^T
\Sigma_0^{-1}(\tilde{\theta} -\theta_0)+O_p
\bigl(n^{-1/2}\bigr),
\end{equation}
and the Lagrange multipliers $\lambda$ at $\theta_0$ can be written as
%
\begin{equation}\label{3-86}
\lambda=\lambda(\theta_0) = \Sigma_0^{-1}(
\tilde{\theta} -\theta_0) +O_p\bigl(n^{-1}
\bigr).
\end{equation}
See \citet{HalLaS90} and \citet{DiCHalRom91N2}. By
(\ref{3-87-2}), (\ref{3-86}) and the Taylor expansion (\ref{3-61}), we have
%
\begin{eqnarray}\label{3-88}
l^*_s(\theta_0)&=& l\bigl(\theta_0+\bigl(
\theta_0'-\theta_0\bigr)\bigr)\nonumber\\
&=&l(\theta
_0)-2n\lambda(\theta_0)^T\bigl(
\theta_0'-\theta_0\bigr) + o_p
\bigl(\bigl\|\bigl(\theta_0'-\theta_0\bigr)\bigr\|
\bigr)
\\
&=& l(\theta_0)-2n\lambda(\theta_0)^T
\frac{b}{2n}(\tilde{\theta}-\theta_0)+O_p
\bigl(n^{-3/2}\bigr).\nonumber
\end{eqnarray}
It follows from (\ref{3-85}), (\ref{3-86}), (\ref{3-88}) and $\tilde
{\theta}-\theta_0=O_p(n^{-1/2})$ that
%
\begin{eqnarray}\label{3-89}
l^*_s(\theta_0) 
&=&l(
\theta_0)-\frac{b}{n} \bigl\{n\bigl[(\tilde{\theta} -
\theta_0)^T\Sigma_0^{-1}
+O_p\bigl(n^{-1}\bigr)\bigr](\tilde{\theta}-
\theta_0) \bigr\}\nonumber\\
&&{}+O_p\bigl(n^{-3/2}\bigr)
\\
&=&l(\theta_0)
\bigl[1-bn^{-1}+O_p\bigl(n^{-3/2}\bigr)\bigr].\nonumber
\end{eqnarray}
This proves (\ref{3-80-1}) and shows that $l^*_s(\theta_0)$ is
equivalent to the BEL in the left-hand side of (\ref{2-72}). Finally,
(\ref{3-81}) follows from (\ref{3-80-1}) and (\ref{2-72}). 
\end{pf*}

\begin{pf*}{Proof of Corollary \ref{cor1}}
It is convenient to view the expansion factor of the first order EEL as
a special case of that for the second order EEL (\ref{3-80}) where
$b=1$ and $\delta(n)=1$. The condition of $\delta(n)=O(n^{-1/2})$
imposed\vspace*{1pt} on the $\delta(n)$ in (\ref{3-80}) is not needed here. Noting
that $l^*(\theta_0)=l(\theta_0)+ O_p(n^{-1})$ and $l(\theta
_0')=l^*(\theta_0)$, equation (\ref{3-83}) in the proof Theorem \ref
{thm3} is now
\[
l\bigl(\theta_0'\bigr)=l(\theta_0)+
O_p\bigl(n^{-1}\bigr).
\]
Thus equation (\ref{3-87-2}) becomes
\[
\theta_0'-\theta_0=\frac{l(\theta_0)}{2n}(
\tilde{\theta}-\theta_0)+O_p\bigl(n^{-5/2}
\bigr).
\]
Using the above equation and following the steps given by (\ref{3-85})
to (\ref{3-89}), we obtain Corollary \ref{cor1}. 
\end{pf*}
\end{appendix}

\section*{Acknowledgements}

Min Tsao is the principal author of this paper. Fan Wu contributed the
simulation results in Section \ref{sec4}.
We thank two anonymous referees and an Associate Editor for their
valuable comments which have led to improvements in this paper. We also
thank the Editor, Professor Peter Hall, and the Associate Editor for
their exceptionally quick handling of this manuscript.



\printaddresses

\end{document}